\documentclass[11pt]{article}

	\newcommand{\mainTitle}{Identity involving symmetric sums of regularized multiple zeta-star values}

	\newcommand{\authorName}{Tomoya Machide}
	\newcommand{\organizationNameFst}{National Institute of Informatics}
	\newcommand{\placeAddressFst}{2-1-2 Hitotsubashi, Chiyoda-ku, Tokyo 101-8430, Japan}
	\newcommand{\emailAddressFst}{machide@nii.ac.jp}
	\newcommand{\organizationNameScd}{JST, ERATO, Kawarabayashi Large Graph Project}
	\newcommand{\departmentNameScd}{c/o Global Research Center for Big Data Mathematics}
	\newcommand{\placeAddressScd}{NII, 2-1-2 Hitotsubashi, Chiyoda-ku, Tokyo 101-8430, Japan}
	
	\newcommand{\MSCname}{11M32 (Primary); 11B73 (Secondary)} 
	\newcommand{\keyWord}{multiple zeta value, multiple zeta-star value, symmetric sum, Bell polynomial} 

\usepackage{graphicx}
\usepackage{geometry}               
\usepackage{ascmac}
\usepackage{amsmath}
\usepackage{amssymb}
\usepackage{amsthm}
\usepackage[colorlinks=true]{hyperref}
\usepackage{arydshln}
\usepackage{multicol}
	
	\DeclareMathOperator*{\OPlus}{\bigoplus}
	\newcommand{\nbk}[3]{#1#3#2}		
	\newcommand{\bgbk}[3]{\bigl{#1}#3\bigr{#2}}	
	\newcommand{\Bgbk}[3]{\Bigl{#1}#3\Bigr{#2}}			
	\newcommand{\bggbk}[3]{\biggl{#1}#3\biggr{#2}}			
	\newcommand{\Bggbk}[3]{\Biggl{#1}#3\Biggr{#2}}
	\newcommand{\autobk}[3]{\left#1#3\right#2}
	\newcommand{\nbkD}[5]{#1#2#5#3#4}		
	\newcommand{\bgbkD}[5]{\bigl{#1}\bigl{#2}#5\bigr{#3}\bigr{#4}}	
	\newcommand{\BgbkD}[5]{\Bigl{#1}\Bigl{#2}#5\Bigr{#3}\Bigr{#4}}	
	\newcommand{\bggbkD}[5]{\biggl{#1}\biggl{#2}#5\biggr{#3}\biggr{#4}}	
	\newcommand{\BggbkD}[5]{\Biggl{#1}\Biggl{#2}#5\Biggr{#3}\Biggr{#4}}	
	\newcommand{\autobkD}[5]{\left#1\left#2#5\right#3\right#4}	
	\newcommand{\mcbk}[4][?]{\ifx n#1\nbk{#2}{#3}{#4}\else\ifx b#1\bgbk{#2}{#3}{#4}\else\ifx B#1\Bgbk{#2}{#3}{#4}\else\ifx g#1\bggbk{#2}{#3}{#4}\else\ifx G#1\Bggbk{#2}{#3}{#4}\else\ifx a#1\autobk{#2}{#3}{#4}\else\ifx !#1{#4}\else#4\fi\fi\fi\fi\fi\fi\fi}
	\newcommand{\mcbkD}[4][?]{\ifx n#1\nbkD{#2}{#2}{#3}{#3}{#4}\else\ifx b#1\bgbkD{#2}{#2}{#3}{#3}{#4}\else\ifx B#1\BgbkD{#2}{#2}{#3}{#3}{#4}\else\ifx g#1\bggbkD{#2}{#2}{#3}{#3}{#4}\else\ifx G#1\BggbkD{#2}{#2}{#3}{#3}{#4}\else\ifx a#1\autobkD{#2}{#2}{#3}{#3}{#4}\else\ifx !#1{#4}\else#4\fi\fi\fi\fi\fi\fi\fi}
	\newcommand{\nsgsb}[1]{#1}		
	\newcommand{\bgsgsb}[1]{\big{#1}}	
	\newcommand{\Bgsgsb}[1]{\Big{#1}}			
	\newcommand{\bggsgsb}[1]{\bigg{#1}}			
	\newcommand{\Bggsgsb}[1]{\Bigg{#1}}
	\newcommand{\mcsgsb}[2][?]{\ifx n#1\nsgsb{#2}\else\ifx b#1\bgsgsb{#2}\else\ifx B#1\Bgsgsb{#2}\else\ifx g#1\bggsgsb{#2}\else\ifx G#1\Bggsgsb{#2}\else#2\fi\fi\fi\fi\fi}
	\newcommand{\myEqSpace}{\,} 	\newlength{\myEqSpaceLen} 	
	\newcommand{\bkR}[2][n]{\mcbk[#1]{(}{)}{#2}}						
	\newcommand{\bkS}[2][n]{\mcbk[#1]{[}{]}{#2}}						
	\newcommand{\bkB}[2][n]{\mcbk[#1]{\{}{\}}{#2}}						
	\newcommand{\bkAll}[4][n]{\mcbk[#1]{#2}{#3}{#4}}
		
	\newcommand{\nFc}[3][n]{#2\bkR[#1]{#3}}					
				
	\newcommand{\idFc}[4][n]{\id{#2}{#3}\bkR[#1]{#4}}			
	\newcommand{\pwFc}[4][n]{\pw{#2}{#3}\bkR[#1]{#4}}			
	\newcommand{\ipFc}[5][n]{\ip{#2}{#3}{#4}\bkR[#1]{#5}}		
		\newcommand{\Fc}{\nFc}		
			
	\newcommand{\Del}{\Delta}
	\newcommand{\gam}{\gamma} 
	\newcommand{\Gam}{\Gamma} 
	
	\newcommand{\ep}{\varepsilon}

	\newcommand{\sig}{\sigma}

	\newcommand{\bfAlp}[1][s]{\ifx s#1{\boldsymbol\alpha}\else{\boldsymbol??none??}\fi}		
	\newcommand{\bfBeta}[1][s]{\ifx s#1{\boldsymbol\beta}\else{\boldsymbol??none??}\fi}		
	\newcommand{\bfDelta}[1][s]{\ifx s#1{\boldsymbol\delta}\else{\boldsymbol??none??}\fi}		
	\newcommand{\bfGam}[1][s]{\ifx s#1{\boldsymbol\gam}\else{\boldsymbol\Gam}\fi}			
	\newcommand{\bfA}[1][s]{\ifx s#1{\bf a}\else{\bf A}\fi}
	\newcommand{\bfB}[1][s]{\ifx s#1{\bf b}\else{\bf B}\fi}
	\newcommand{\bfE}[1][s]{\ifx s#1{\bf e}\else{\bf E}\fi}
	\newcommand{\bfH}[1][s]{\ifx s#1{\bf h}\else{\bf H}\fi}
	\newcommand{\bfI}[1][s]{\ifx s#1{\bf i}\else{\bf I}\fi}
	\newcommand{\bfJ}[1][s]{\ifx s#1{\bf j}\else{\bf J}\fi}
	\newcommand{\bfK}[1][s]{\ifx s#1{\bf k}\else{\bf K}\fi}
	\newcommand{\bfL}[1][s]{\ifx s#1{\bf l}\else{\bf L}\fi}
	\newcommand{\bfM}[1][s]{\ifx s#1{\bf m}\else{\bf M}\fi}
	\newcommand{\bfN}[1][s]{\ifx s#1{\bf n}\else{\bf N}\fi}
	\newcommand{\bfS}[1][s]{\ifx s#1{\bf s}\else{\bf S}\fi}
	\newcommand{\bfU}[1][s]{\ifx s#1{\bf u}\else{\bf U}\fi}
	\newcommand{\bfV}[1][s]{\ifx s#1{\bf v}\else{\bf V}\fi}
	\newcommand{\bfW}[1][s]{\ifx s#1{\bf w}\else{\bf W}\fi}
	\newcommand{\bfX}[1][s]{\ifx s#1{\bf x}\else{\bf X}\fi}
	\newcommand{\bfY}[1][s]{\ifx s#1{\bf y}\else{\bf Y}\fi}
	\newcommand{\bfZ}[1][s]{\ifx s#1{\bf z}\else{\bf Z}\fi}
	\newcommand{\bfEp}[1][s]{\ifx s#1{\boldsymbol\ep}\else{\boldsymbol{\mathcal{E}}}\fi}
	
	\newcommand{\mo}{(-1)}
	\newcommand{\mVert}[1][n]{{\,\mcsgsb[#1]{\vert}\,}}		
	
	\newcommand{\SetO}[2][n]{\bkB[#1]{#2}}
	\newcommand{\SetT}[3][n]{\bkB[#1]{#2\mVert#3}}
		\newcommand{\Set}{\SetO}

	\newcommand{\setQ}{\mathbb{Q}}	
	\newcommand{\setR}{\mathbb{R}}	
	
	\newcommand{\setF}[1][?]{\ifx ?#1\mathbb{F}\else\mathbb{F}_{#1}\fi}
	
	\newcommand{\matI}[1][?]{\ifx #1?I\else I_{#1}\fi}	
	
	\newcommand{\gpSym}[2][?]{S_{#2}}	
	\newcommand{\gpKleinF}[1][?]{V}
		
	\newcommand{\gpu}[1][?]{\ifx?#1e\else e_{#1}\fi}		

	\newcommand{\vPack}[1][10]{\vspace{-#1pt}}
	
	\newcommand{\lnA}[1][]{&  &}
	\newcommand{\lnP}[1]{\myEqSpace#1\myEqSpace}
	\newcommand{\lnAP}[2][]{& #2 &}
	
	\newcommand{\lnAHP}[2][\nonumber]{#1 \\ & #2 &}

		\newcommand{\slnAH}[1][?]{\\}
		
	\newcommand{\refEq}[1]{(\ref{#1})}	
	\newcommand{\refEqA}[1]{(#1)}	
	\newcommand{\pcstSpForRefThm}{\;}		
	\newcommand{\refHL}[2]{#1\pcstSpForRefThm\ref{#2}}		
	\newcommand{\refHLm}[3][?]{\ifx?#1#2\pcstSpForRefThm#3\else#2#3\fi}
	\newcommand{\refThm}[2][?]{\ifx?#1\refHL{Theorem}{#2}\else\ifx s#1\refHL{Theorems}{#2}\else{[argument error]}\fi\fi}
	\newcommand{\refProp}[2][?]{\ifx?#1\refHL{Proposition}{#2}\else\ifx s#1\refHL{Propositions}{#2}\else{[argument error]}\fi\fi}
	\newcommand{\refLem}[2][?]{\ifx?#1\refHL{Lemma}{#2}\else\ifx s#1\refHL{Lemmas}{#2}\else{[argument error]}\fi\fi}
	\newcommand{\refCor}[2][?]{\ifx?#1\refHL{Corollary}{#2}\else\ifx s#1\refHL{Corollaries}{#2}\else{[argument error]}\fi\fi}
	\newcommand{\refDef}[2][?]{\ifx?#1\refHL{Definition}{#2}\else\ifx s#1\refHL{Definitions}{#2}\else{[argument error]}\fi\fi}
	\newcommand{\refRem}[2][?]{\ifx?#1\refHL{Remark}{#2}\else\ifx s#1\refHL{Remarks}{#2}\else{[argument error]}\fi\fi}
	\newcommand{\refTab}[2][?]{\ifx?#1\refHL{Table}{#2}\else\ifx s#1\refHL{Tables}{#2}\else{[argument error]}\fi\fi}
	\newcommand{\refSec}[2][?]{\ifx?#1\refHL{Section}{#2}\else\ifx s#1\refHL{Sections}{#2}\else{[argument error]}\fi\fi}
	\newcommand{\refApp}[2][?]{\ifx?#1\refHL{Appendix}{#2}\else\ifx s#1\refHL{Appendices}{#2}\else{[argument error]}\fi\fi}
		
	\newcommand{\refThmA}[2][?]{\ifx?#1\refHLm{Theorem}{#2}\else \refHLm{Theorems}{#2}\fi}
	\newcommand{\refPropA}[2][?]{\ifx?#1 \refHLm[#1]{Proposition}{#2}\else \refHLm{Propositions}{#2}\fi}
	\newcommand{\refLemA}[2][?]{\ifx?#1\refHLm{Lemma}{#2}\else \refHLm{Lemmas}{#2}\fi}
	\newcommand{\refCorA}[2][?]{\ifx?#1\refHLm{Corollary}{#2}\else \refHLm{Corollaries}{#2}\fi}
	\newcommand{\refDefA}[2][?]{\ifx?#1\refHLm{Definition}{#2}\else \refHLm{Definitionss}{#2}\fi}
	\newcommand{\refRemA}[2][?]{\ifx?#1\refHLm{Remark}{#2}\else \refHLm{Remarks}{#2}\fi}
	\newcommand{\refSecA}[2][?]{\ifx?#1\refHLm{Section}{#2}\else \refHLm{Sections}{#2}\fi}
	\newcommand{\refAppA}[2][?]{\ifx?#1\refHLm{Appendix}{#2}\else \refHLm{Appendices}{#2}\fi}

	\newcommand{\frc}[3][?]{\ifx s#1#3/#2\else\ifx b#1(#3)/#2\else\ifx d#1\dfrac{#3}{#2}\else\ifx t#1\tfrac{#3}{#2}\else\frac{#3}{#2}\fi\fi\fi\fi}
	\newcommand{\bnm}[3][?]{\binom{#3}{#2}}

	\newcommand{\pw}[3][?]{\ifx!#3{#2}^{#3}\else#2^{#3}\fi}
	\newcommand{\id}[3][?]{#2_{#3}}
	\newcommand{\ip}[4][?]{{#2}_{#3}^{#4}}
	\newcommand{\pwR}[3][a]{\ifx!#1{\bkR[#1]{#2}}^{#3}\else\bkR[#1]{#2}^{#3}\fi}
	\newcommand{\pwB}[3][a]{\ifx!#1{\bkB[#1]{#2}}^{#3}\else\bkB[#1]{#2}^{#3}\fi}
	\newcommand{\pwS}[3][a]{\ifx!#1{\bkS[#1]{#2}}^{#3}\else\bkS[#1]{#2}^{#3}\fi}

	\newcommand{\tpT}[3][a]{ {#2}\atop \bkR[#1]{#3} }

	\newcommand{\nSmO}[2][?]{\ifx l#1\sum\limits_{#2}\else\ifx t#1{\textstyle\sum\limits_{#2}}\else\sum_{#2}\fi\fi}
	\newcommand{\nSmT}[3][?]{\ifx l#1\sum\limits_{#2}^{#3}\else\if t#1{\textstyle\sum\limits_{#2}^{#3}}\else\sum_{#2}^{#3}\fi\fi}	
	\newcommand{\nSmN}[1][?]{\ifx l#1\sum\limits\else\ifx t#1{\textstyle\sum\limits}\else\sum\fi\fi}
	\newcommand{\pSm}[2][?]{\ifx t#1 \sum_{#2}^{\prime} \else \sideset{}{^\prime}\sum_{#2} \fi}
	\newcommand{\pSmT}[3][?]{\ifx t#1 \sum_{#2}^{\prime#3} \else \sideset{}{^\prime}\sum_{#2}^{#3} \fi}	
	\newcommand{\pSmN}[1][?]{\ifx t#1 \sum^{\prime} \else \sideset{}{^\prime}\sum \fi}
	\newcommand{\dSm}[2][?]{\ifx t#1 \sum_{#2}^{\dagger} \else \sideset{}{^\dagger}\sum_{#2} \fi}
	\newcommand{\dSmT}[3][?]{\ifx t#1 \sum_{#2}^{\dagger#3} \else \sideset{}{^\dagger}\sum_{#2}^{#3} \fi}	
	\newcommand{\dSmN}[1][?]{\ifx t#1 \sum^{\dagger} \else \sideset{}{^\dagger}\sum \fi}
	\newcommand{\tpTSm}[3][?]{\nSmO[#1]{\tpT{#2}{#3}}}

		\newcommand{\Sm}{\nSmO}			\newcommand{\SmT}{\nSmT}			
		\newcommand{\tpSm}{\tpTSm}
		
	\newcommand{\nPd}[2][?]{\ifx l#1 \prod\limits_{#2} \else \prod_{#2} \fi}
	\newcommand{\nPdT}[3][?]{\ifx l#1 \prod\limits_{#2}^{#3} \else \prod_{#2}^{#3} \fi}

		\newcommand{\PdT}{\nPdT}

	\newcommand{\nOPs}[2][?]{\ifx l#1 \OPlus\limits_{#2} \else \OPlus_{#2} \fi}
	\newcommand{\nOPsT}[3][?]{\ifx l#1 \OPlus\limits_{#2}^{#3} \else \OPlus_{#2}^{#3} \fi}	
	\newcommand{\pOPs}[2][?]{\ifx t#1 \OPlus_{#2}^{\prime} \else \sideset{}{^\prime}\OPlus_{#2} \fi}
	\newcommand{\pOPsT}[3][?]{\ifx t#1 \OPlus_{#2}^{\prime#3} \else \sideset{}{^\prime}\OPlus_{#2}^{#3} \fi}

	\newcommand{\nIs}[2][?]{\ifx l#1 \bigcap\limits_{#2}\else\ifx b#1 \bigcap_{#2}\else{\textstyle\bigcap\limits_{#2}}\fi\fi}
	\newcommand{\nIsT}[3][?]{\ifx l#1 \bigcap\limits_{#2}^{#3}\else\ifx b#1 \bigcap_{#2}^{#3}\else{\textstyle\bigcap\limits_{#2}^{#3}}\fi\fi}	
	\newcommand{\pIs}[2][?]{\ifx t#1 \bigcap_{#2}^{\prime} \else \sideset{}{^\prime}\bigcap_{#2} \fi}
	\newcommand{\pIsT}[3][?]{\ifx t#1 \bigcap_{#2}^{\prime#3} \else \sideset{}{^\prime}\bigcap_{#2}^{#3} \fi}

	\newcommand{\nUn}[2][?]{\ifx L#1 \bigcup\limits_{#2}\else\ifx b#1 \bigcup_{#2}\else\ifx t#1{\textstyle\bigcup_{#2}}\else{\textstyle\bigcup_{#2}}\fi\fi\fi}
	\newcommand{\nUnT}[3][?]{\ifx L#1 \bigcup\limits_{#2}^{#3}\else\ifx b#1 \bigcup_{#2}^{#3}\else\ifx t#1{\textstyle\bigcup_{#2}^{#3}}\else{\textstyle\bigcup\limits_{#2}^{#3}}\fi\fi\fi}	
	\newcommand{\pUn}[2][?]{\ifx t#1 \bigcup_{#2}^{\prime} \else \sideset{}{^\prime}\bigcup_{#2} \fi}
	\newcommand{\pUnT}[3][?]{\ifx t#1 \bigcup_{#2}^{\prime#3} \else \sideset{}{^\prime}\bigcup_{#2}^{#3} \fi}	
	\newcommand{\dUn}[2][?]{\ifx L#1 \bigsqcup\limits_{#2}\else\ifx b#1 \bigsqcup_{#2}\else\ifx t#1{\textstyle\bigsqcup_{#2}}\else{\textstyle\bigsqcup\limits_{#2}}\fi\fi\fi}
	\newcommand{\dUnT}[3][?]{\ifx L#1 \bigsqcup\limits_{#2}^{#3}\else\ifx b#1 \bigsqcup_{#2}^{#3}\else\ifx t#1{\textstyle\bigsqcup_{#2}^{#3}}\else{\textstyle\bigsqcup\limits_{#2}^{#3}}\fi\fi\fi}

	\newcommand{\nLm}[2][?]{\ifx l#1 \lim\limits_{#2} \else \lim_{#2} \fi}
	\newcommand{\iLm}[2][?]{\ifx l#1 \liminf\limits_{#2} \else \liminf_{#2} \fi}
	\newcommand{\sLm}[2][?]{\ifx l#1 \limsup\limits_{#2} \else \limsup_{#2} \fi}
		
	\newcommand{\glcondEnvLineHead}[1]{ \ifx*#1 \begin{eqnarray*} \else \begin{eqnarray}  \label{#1} \fi }
	\newcommand{\glcondEnvLineTail}[1]{ \ifx*#1 \end{eqnarray*} \else \end{eqnarray} \fi }
	\newcommand{\glcondDis}[1]{\ifx d#1 \displaystyle \fi}
	\newcommand{\glcmdEqShift}{\hspace{-20pt}}
	\newcommand{\glcmdHLineCWiden}{\rule{0cm}{15pt}}	\newcommand{\glcdH}{\glcmdHLineCWiden}
	\newcommand{\lccondPar}[1]{\ifx#1p \\ \fi}

		\newcommand{\envMO}[2][*]{$\ifx d#1 \displaystyle \fi#2$}
		\newcommand{\envMT}[3][*]{$\ifx d#1 \displaystyle \fi#2=#3$}
		\newcommand{\envMTDef}[3][*]{$\ifx d#1 \displaystyle \fi#2:=#3$}
		\newcommand{\envMTPt}[4][*]{$\ifx d#1 \displaystyle \fi#3#2#4$}
			\newcommand{\envM}{\envMT}

		\newcommand{\envMTh}[4][*]{$\ifx d#1 \displaystyle \fi#2=#3=#4$}
		\newcommand{\envMThDef}[4][*]{$\ifx d#1 \displaystyle \fi#2:=#3=#4$}
		\newcommand{\envMThPt}[5][*]{$\ifx d#1 \displaystyle \fi#3#2#4#2#5$}
		\newcommand{\envMThPte}[6][*]{$\ifx d#1 \displaystyle \fi#2#3#4#5#6$}
		\newcommand{\envMF}[5][*]{$\ifx d#1 \displaystyle \fi#2=#3=#4=#5$}
		\newcommand{\envMFPt}[6][*]{$\ifx d#1 \displaystyle \fi#3#2#4#2#5#2#6$}
		
		\newcommand{\envHLineT}[3][*]{ \glcondEnvLineHead{#1} #2&=&#3\glcondEnvLineTail{#1} }
		\newcommand{\envHLineTDef}[3][*]{ \glcondEnvLineHead{#1} #2&:=&#3\glcondEnvLineTail{#1} }

			\newcommand{\envHLine}{\envHLineT}
			\newcommand{\envHLineDef}{\envHLineTDef}

		\newcommand{\envPLine}[2][*]{\glcondEnvLineHead{#1} #2\glcondEnvLineTail{#1}}
		
		\newcommand{\envOTLine}[4][*]{\glcondEnvLineHead{#1} #2\lnAP{=}#3\lnP{=}#4\glcondEnvLineTail{#1}}

			\newcommand{\envOTLineTh}{\envOTLine}
			
	\newcommand{\lcparaCase}{\vspace{3pt}}

	\newcommand{\matu}[1][?]{\ifx#1?I\else I_{#1}\fi}
	\newcommand{\vlA}[2][n]{\bkAll[#1]{|}{|}{#2}}
				
	\newcommand{\nmSet}[2][n]{\vlA[#1]{#2}}	
		\newcommand{\vlNS}{\nmSet}

	\newcommand{\exx}[2][n]{ \Fc[#1]{\exp}{#2} }
	
	\newcommand{\exN}[1]{ e^{#1}}
	
	\newcommand{\fcGam}[2][n]{\Fc[#1]{\Gam}{#2}}

		\newcommand{\fcG}{\fcGam}

	\newcommand{\sTx}[2][?]{ \ifx t#1{\tiny #2} \else \ifx s#1{\scriptsize #2} \else \ifx f#1{\footnotesize #2} \else \ifx S#1{\small #2} \else \ifx n#1{\normalsize #2} \else \ifx l#1{\large #2} \else \ifx L#1{\Large #2} \else \ifx R#1{\LARGE #2} \else \ifx h#1{\huge #2} \else \ifx H#1{\Huge #2} \else \ifx ?#1 #2 \else #2 \fi\fi\fi\fi\fi\fi\fi\fi\fi\fi\fi }
	\newcommand{\bfTx}[1]{{\bf#1}}
	
	\newcommand{\osMTx}[3][?]{\overset{#3}{#2}}
	\newcommand{\usbMTx}[3][?]{\underset{#2}{\underbrace{#3}}}

		\newcommand{\osTx}{\osMTx}
		\newcommand{\usbTx}{\usbMTx}
		
	\newcommand{\raTx}[3][?]{\raisebox{#2pt}[0pt][0pt]{\ifx d#1\displaystyle\fi#3}}
	\newcommand{\raMTx}[3][?]{\raisebox{#2pt}[0pt][0pt]{$\ifx d#1\displaystyle\fi#3$}}

	\newcommand{\envCenter}[2][*]{\ifx*#1\begin{center}\else\begin{center}[#1]\fi #2\end{center}}
	\newcommand{\envFlushleft}[2][*]{\ifx*#1\begin{flushleft}\else\begin{flushleft}[#1]\fi #2\end{flushleft}}
	\newcommand{\envFlushright}[2][*]{\ifx*#1\begin{flushright}\else\begin{flushright}[#1]\fi #2\end{flushright}}
		
	\newcommand{\envItemIm}[2][*]{\ifx*#1\begin{itemize}\else\begin{itemize}[#1]\fi #2\end{itemize}}
	\newcommand{\envItemDp}[2][*]{\ifx*#1\begin{description}\else\begin{description}[#1]\fi #2\end{description}}
	\newcommand{\envItemEm}[2][*]{\ifx*#1\begin{enumerate}\else\begin{enfumerate}[#1]\fi #2\end{enumerate}}

	\newcommand{\envMultCol}[3][*]{\ifx1#2#3\else\begin{multicols}{#2}\ifx*#1\else\mbox{}\vspace{-#1pt}\fi#3\end{multicols}\fi}
	
\theoremstyle{plain}
\newtheorem{theorem}{THEOREM}[section]
\newtheorem{proposition}[theorem]{PROPOSITION}
\newtheorem{lemma}[theorem]{LEMMA}
\newtheorem{corollary}[theorem]{COROLLARY}
\theoremstyle{definition}

\newtheorem{example}[theorem]{EXAMPLE}
\theoremstyle{remark}
\newtheorem{remark}[theorem]{REMARK}
\theoremstyle{plain}

\theoremstyle{definition}

\theoremstyle{remark}

\theoremstyle{plain}

\theoremstyle{definition}

\theoremstyle{remark}

\allowdisplaybreaks[4]
\numberwithin{equation}{section}

	\newcommand{\lccondBibitem}[3][]{ \if ?#2 \bibitem{#3} \else \bibitem[#2]{#3} \fi}
	\newcommand{\refPaper}[8][?]{
			\lccondBibitem{#1}{#2}
				#3,			
				\emph{#4}, 	
				#5\ 			
				{\bf #6}		
				(#7),			
				#8.			
		}
	
	\newcommand{\refPaperRep}[9][?]{
			\lccondBibitem{#1}{#2}
				#3,			
				\emph{#4}, 	
				#5\ 			
				{\bf #6}		
				(#7),			
				#8			
				; reprinted in #9	
		}

	\newcommand{\refPaperAlm}[5][?]{
			\lccondBibitem{#1}{#2}
				#3,	 		
				\emph{#4}, 	
				#5		
		}

	\newcommand{\glcondEnvLineTailPd}[1]{.\ifx*#1 \end{eqnarray*} \else \end{eqnarray} \fi  }
	\newcommand{\glcondEnvLineTailCm}[1]{,\ifx*#1 \end{eqnarray*} \else \end{eqnarray} \fi }
	\newcommand{\prcondEnvEqSpHead}[1]{ \ifx*#1 \begin{equation*}[ERROR] \else \begin{equation}  \label{#1} \fi  }
	\newcommand{\prcondEnvEqSpTail}[1]{\ifx*#1 [ERROR]\end{equation*} \else \end{equation} \fi }

	\newcommand{\envProof}[2][?]{ \par\mbox{}\vspace{-5pt}\\ \ifx?#1\emph{Proof.}\else\emph{Proof of #1.}\fi \ #2 \hfill $\Box$\\ \par}
	
		\newcommand{\envLineTCm}[3][*]{ \glcondEnvLineHead{#1} & &\glcmdEqShift#2\nonumber\\&=&#3 \glcondEnvLineTailCm{#1} }

			\newcommand{\envLineCm}{\envLineTCm}
			
		\newcommand{\envLineTCmPt}[4][*]{\glcondEnvLineHead{#1} & &\glcmdEqShift#3\nonumber\\&#2&#4 \glcondEnvLineTailCm{#1}}

			\newcommand{\envLineCmPt}{\envLineTCmPt}

		\newcommand{\envLineThPd}[4][*]{ \glcondEnvLineHead{#1} & &\glcmdEqShift#2\nonumber\\&=&#3\nonumber \\&=&#4 \glcondEnvLineTailPd{#1} }
		\newcommand{\envLineThCm}[4][*]{ \glcondEnvLineHead{#1} & &\glcmdEqShift#2\nonumber\\&=&#3\nonumber \\&=&#4 \glcondEnvLineTailCm{#1} }

		\newcommand{\envLineFPd}[5][*]{ \glcondEnvLineHead{#1} & &\glcmdEqShift#2\nonumber\\&=&#3\nonumber \\&=&#4\nonumber \\&=&#5 \glcondEnvLineTailPd{#1} }
		
		\newcommand{\envHLineTPd}[3][*]{ \glcondEnvLineHead{#1} #2&=&#3\glcondEnvLineTailPd{#1} }
		\newcommand{\envHLineTDefPd}[3][*]{ \glcondEnvLineHead{#1} #2&:=&#3\glcondEnvLineTailPd{#1} }
		\newcommand{\envHLineTCm}[3][*]{ \glcondEnvLineHead{#1} #2&=&#3\glcondEnvLineTailCm{#1} }
		\newcommand{\envHLineTCmDef}[3][*]{ \glcondEnvLineHead{#1} #2&:=&#3\glcondEnvLineTailCm{#1} }
		\newcommand{\envHLineTCmPt}[4][*]{\glcondEnvLineHead{#1} #3&#2&#4\glcondEnvLineTailCm{#1}}

			\newcommand{\envHLinePd}{\envHLineTPd}
			\newcommand{\envHLineDefPd}{\envHLineTDefPd}
			\newcommand{\envHLineCm}{\envHLineTCm}
			\newcommand{\envHLineCmDef}{\envHLineTCmDef}
			\newcommand{\envHLineCmPt}{\envHLineTCmPt}

		\newcommand{\envHLineFPd}[5][*]{ \glcondEnvLineHead{#1} #2&=&#3\nonumber\\&=&#4\nonumber \\&=&#5 \glcondEnvLineTailPd{#1} }

		\newcommand{\envHLineFCmPte}[8][*]{\glcondEnvLineHead{#1} #2&#3&#4\nonumber\\&#5&#6\nonumber \\&#7&#8\glcondEnvLineTailCm{#1}}
		
		\newcommand{\envHLineCFCmNme}[5][*]{\begin{eqnarray} #2&=&#3,\\\glcdH#4&=&#5 \glcondEnvLineTailCm{?} }
		\newcommand{\envHLineCFNmePd}[5][*]{\begin{eqnarray} #2&=&#3,\\\glcdH#4&=&#5 \glcondEnvLineTailPd{?} }
		\newcommand{\envHLineCFCmDefNme}[5][*]{\begin{eqnarray} #2&:=&#3,\\\glcdH#4&:=&#5 \glcondEnvLineTailCm{?} }
		\newcommand{\envHLineCFDefNmePd}[5][*]{\begin{eqnarray} #2&:=&#3,\\\glcdH#4&:=&#5 \glcondEnvLineTailPd{?} }

		\newcommand{\envHLineCFNmePdPt}[6][*]{\begin{eqnarray}#3&#2&#4,\\\glcdH#5&#2&#6\glcondEnvLineTailPd{?}}
		\newcommand{\envHLineCFCmNmePt}[6][*]{\begin{eqnarray}#3&#2&#4,\\\glcdH#5&#2&#6\glcondEnvLineTailCm{?}}
		\newcommand{\envHLineCFNmePdPte}[7][*]{\begin{eqnarray}#2&#3&#4,\\\glcdH#5&#6&#7\glcondEnvLineTailPd{?}}
		\newcommand{\envHLineCFCmNmePte}[7][*]{\begin{eqnarray}#2&#3&#4,\\\glcdH#5&#6&#7\glcondEnvLineTailCm{?}}

		\newcommand{\envHLineCFLaaCm}[5][*]{\envPLineCm[#1]{#2\lnP{=}#3\qquad\text{and}\qquad#4\lnP{=}#5}}

		\newcommand{\envHLineCSNmePd}[7][*]{\begin{eqnarray} #2&=&#3,\\\glcdH#4&=&#5,\\\glcdH#6&=&#7\glcondEnvLineTailPd{?}}
		\newcommand{\envHLineCSDefNmePd}[7][*]{\begin{eqnarray} #2&:=&#3,\\\glcdH#4&:=&#5,\\\glcdH#6&:=&#7\glcondEnvLineTailPd{?}}
		\newcommand{\envHLineCSCmNme}[7][*]{\begin{eqnarray} #2&=&#3,\\\glcdH#4&=&#5,\\\glcdH#6&=&#7\glcondEnvLineTailCm{?}}
		\newcommand{\envHLineCSCmDefNme}[7][*]{\begin{eqnarray} #2&:=&#3,\\\glcdH#4&:=&#5,\\\glcdH#6&:=&#7\glcondEnvLineTailCm{?}}

		\newcommand{\envHLineCSNmePdPt}[8][*]{\begin{eqnarray}#3&#2&#4,\\\glcdH#5&#2&#6,\\\glcdH#7&#2&#8\glcondEnvLineTailPd{?}}
		\newcommand{\envHLineCSCmNmePt}[8][*]{\begin{eqnarray}#3&#2&#4,\\\glcdH#5&#2&#6,\\\glcdH#7&#2&#8\glcondEnvLineTailCm{?}}
		\newcommand{\envHLineCSNmePdPte}[9][*]{\begin{eqnarray}#2&#3&#4,\\\glcdH#5&#6&#7,\\\glcdH#8&#2&#9\glcondEnvLineTailPd{?}}
		\newcommand{\envHLineCSCmNmePte}[9][*]{\begin{eqnarray}#2&#3&#4,\\\glcdH#5&#6&#7,\\\glcdH#8&#2&#9\glcondEnvLineTailCm{?}}
		
		\newcommand{\envHLineCENmePd}[9][*]{\begin{eqnarray} #2&=&#3,\\\glcdH#4&=&#5,\\\glcdH#6&=&#7,\\\glcdH#8&=&#9\glcondEnvLineTailPd{?}}
		\newcommand{\envHLineCEDefNmePd}[9][*]{\begin{eqnarray} #2&:=&#3,\\\glcdH#4&:=&#5,\\\glcdH#6&:=&#7,\\\glcdH#8&:=&#9\glcondEnvLineTailPd{?}}
		\newcommand{\envHLineCECmNme}[9][*]{\begin{eqnarray} #2&=&#3,\\\glcdH#4&=&#5,\\\glcdH#6&=&#7,\\\glcdH#8&=&#9\glcondEnvLineTailCm{?}}
		\newcommand{\envHLineCECmDefNme}[9][*]{\begin{eqnarray} #2&:=&#3,\\\glcdH#4&:=&#5,\\\glcdH#6&:=&#7,\\\glcdH#8&:=&#9\glcondEnvLineTailCm{?}}

			\newcommand{\pccondPaForPar}[1]{\ifx#1p \\\glcdH \fi}
			\newcommand{\pccondPaForNonnum}[1]{\ifx#1* \nonumber \fi}

		\newcommand{\envPLineCm}[2][*]{\glcondEnvLineHead{#1} #2\glcondEnvLineTailCm{#1}}
	
		\newcommand{\envOTLinePd}[4][*]{\glcondEnvLineHead{#1} #2\lnAP{=}#3\lnP{=}#4.\glcondEnvLineTail{#1}}

			\newcommand{\envOTLineThPd}{\envOTLinePd}

		\newcommand{\envMOCm}[2][*]{$\ifx d#1 \displaystyle \fi#2$,}
		\newcommand{\envMOPd}[2][*]{$\ifx d#1 \displaystyle \fi#2$.}
		
		\newcommand{\envMTCm}[3][*]{$\ifx d#1 \displaystyle \fi#2=#3$,}
		\newcommand{\envMTPd}[3][*]{$\ifx d#1 \displaystyle \fi#2=#3$.}
		\newcommand{\envMTCmDef}[3][*]{$\ifx d#1 \displaystyle \fi#2:=#3$,}
		\newcommand{\envMTDefPd}[3][*]{$\ifx d#1 \displaystyle \fi#2:=#3$.}
		\newcommand{\envMTCmPt}[4][*]{$\ifx d#1 \displaystyle \fi#3#2#4$,}
		\newcommand{\envMTPdPt}[4][*]{$\ifx d#1 \displaystyle \fi#3#2#4$.}
			
			\newcommand{\envMPd}{\envMTPd}

		\newcommand{\envMThCm}[4][*]{$\ifx d#1 \displaystyle \fi#2=#3=#4$,}
		\newcommand{\envMThPd}[4][*]{$\ifx d#1 \displaystyle \fi#2=#3=#4$.}
		\newcommand{\envMThCmPt}[5][*]{$\ifx d#1 \displaystyle \fi#3#2#4#2#5$,}
		\newcommand{\envMThPdPt}[5][*]{$\ifx d#1 \displaystyle \fi#3#2#4#2#5$.}
		\newcommand{\envMFCm}[5][*]{$\ifx d#1 \displaystyle \fi#2=#3=#4=#5$,}
		\newcommand{\envMFPd}[5][*]{$\ifx d#1 \displaystyle \fi#2=#3=#4=#5$.}
		\newcommand{\envMFCmPt}[6][*]{$\ifx d#1 \displaystyle \fi#3#2#4#2#5#2#6$,}
		\newcommand{\envMFPdPt}[6][*]{$\ifx d#1 \displaystyle \fi#3#2#4#2#5#2#6$.}

		\newcommand{\prcondHLCPNm}{\hspace{-1pt}}
		\newcommand{\envHLineCFCmNm}[5][*]{ \begin{equation}\begin{split} \ifx*#1 \text{[ERROR;need label name]} \else \label{#1} \fi #2&\prcondHLCPNm\lnP{=}\prcondHLCPNm#3,\\#4&\prcondHLCPNm\lnP{=}\prcondHLCPNm#5, \end{split}\end{equation} }
		\newcommand{\envHLineCFNm}[5][*]{ \begin{equation}\begin{split} \ifx*#1 \text{[ERROR;need label name]} \else \label{#1} \fi #2&\prcondHLCPNm\lnP{=}\prcondHLCPNm#3\\#4&\prcondHLCPNm\lnP{=}\prcondHLCPNm#5, \end{split}\end{equation} }
		\newcommand{\envHLineCFNmPd}[5][*]{ \begin{equation}\begin{split} \ifx*#1 \text{[ERROR;need label name]} \else \label{#1} \fi #2&\prcondHLCPNm\lnP{=}\prcondHLCPNm#3,\\#4&\prcondHLCPNm\lnP{=}\prcondHLCPNm#5. \end{split}\end{equation} }
		\newcommand{\envHLineCFCmDefNm}[5][*]{ \begin{equation}\begin{split} \ifx*#1 \text{[ERROR;need label name]} \else \label{#1} \fi #2&\prcondHLCPNm\lnP{:=}\prcondHLCPNm#3,\\#4&\prcondHLCPNm\lnP{:=}\prcondHLCPNm#5, \end{split}\end{equation} }
		\newcommand{\envHLineCFDefNm}[5][*]{ \begin{equation}\begin{split} \ifx*#1 \text{[ERROR;need label name]} \else \label{#1} \fi #2&\prcondHLCPNm\lnP{:=}\prcondHLCPNm#3\\#4&\prcondHLCPNm\lnP{:=}\prcondHLCPNm#5, \end{split}\end{equation} }
		\newcommand{\envHLineCFDefNmPd}[5][*]{ \begin{equation}\begin{split} \ifx*#1 \text{[ERROR;need label name]} \else \label{#1} \fi #2&\prcondHLCPNm\lnP{:=}\prcondHLCPNm#3,\\#4&\prcondHLCPNm\lnP{:=}\prcondHLCPNm#5. \end{split}\end{equation} }
		\newcommand{\envHLineCSCmNm}[7][*]{ \begin{equation}\begin{split} \ifx*#1 \text{[ERROR;need label name]} \else \label{#1} \fi #2&\prcondHLCPNm\lnP{=}\prcondHLCPNm#3,\\#4&\prcondHLCPNm\lnP{=}\prcondHLCPNm#5,\\#6&\prcondHLCPNm\lnP{=}\prcondHLCPNm#7 \end{split}\end{equation} }
		\newcommand{\envHLineCSNm}[7][*]{ \begin{equation}\begin{split} \ifx*#1 \text{[ERROR;need label name]} \else \label{#1} \fi #2&\prcondHLCPNm\lnP{=}\prcondHLCPNm#3\\#4&\prcondHLCPNm\lnP{=}\prcondHLCPNm#5\\#6&\prcondHLCPNm\lnP{=}\prcondHLCPNm#7 \end{split}\end{equation} }
		\newcommand{\envHLineCSNmPd}[7][*]{ \begin{equation}\begin{split} \ifx*#1 \text{[ERROR;need label name]} \else \label{#1} \fi #2&\prcondHLCPNm\lnP{=}\prcondHLCPNm#3,\\#4&\prcondHLCPNm\lnP{=}\prcondHLCPNm#5,\\#6&\prcondHLCPNm\lnP{=}\prcondHLCPNm#7. \end{split}\end{equation} }
		\newcommand{\envHLineCSCmDefNm}[7][*]{ \begin{equation}\begin{split} \ifx*#1 \text{[ERROR;need label name]} \else \label{#1} \fi #2&\prcondHLCPNm\lnP{:=}\prcondHLCPNm#3,\\#4&\prcondHLCPNm\lnP{:=}\prcondHLCPNm#5,\\#6&\prcondHLCPNm\lnP{:=}\prcondHLCPNm#7 \end{split}\end{equation} }
		\newcommand{\envHLineCSDefNm}[7][*]{ \begin{equation}\begin{split} \ifx*#1 \text{[ERROR;need label name]} \else \label{#1} \fi #2&\prcondHLCPNm\lnP{:=}\prcondHLCPNm#3\\#4&\prcondHLCPNm\lnP{:=}\prcondHLCPNm#5\\#6&\prcondHLCPNm\lnP{:=}\prcondHLCPNm#7 \end{split}\end{equation} }
		\newcommand{\envHLineCSDefNmPd}[7][*]{ \begin{equation}\begin{split} \ifx*#1 \text{[ERROR;need label name]} \else \label{#1} \fi #2&\prcondHLCPNm\lnP{:=}\prcondHLCPNm#3,\\#4&\prcondHLCPNm\lnP{:=}\prcondHLCPNm#5,\\#6&\prcondHLCPNm\lnP{:=}\prcondHLCPNm#7. \end{split}\end{equation} }
		\newcommand{\envHLineCECmNm}[9][*]{ \begin{equation}\begin{split} \ifx*#1 \text{[ERROR;need label name]} \else \label{#1} \fi #2&\prcondHLCPNm\lnP{=}\prcondHLCPNm#3,\\#4&\prcondHLCPNm\lnP{=}\prcondHLCPNm#5,\\#6&\prcondHLCPNm\lnP{=}\prcondHLCPNm#7,\\#8&\prcondHLCPNm\lnP{=}\prcondHLCPNm#9,  \end{split}\end{equation} }
		\newcommand{\envHLineCENm}[9][*]{ \begin{equation}\begin{split} \ifx*#1 \text{[ERROR;need label name]} \else \label{#1} \fi #2&\prcondHLCPNm\lnP{=}\prcondHLCPNm#3\\#4&\prcondHLCPNm\lnP{=}\prcondHLCPNm#5\\#6&\prcondHLCPNm\lnP{=}\prcondHLCPNm#7\\#8&\prcondHLCPNm\lnP{=}\prcondHLCPNm#9  \end{split}\end{equation} }
		\newcommand{\envHLineCENmPd}[9][*]{ \begin{equation}\begin{split} \ifx*#1 \text{[ERROR;need label name]} \else \label{#1} \fi #2&\prcondHLCPNm\lnP{=}\prcondHLCPNm#3,\\#4&\prcondHLCPNm\lnP{=}\prcondHLCPNm#5,\\#6&\prcondHLCPNm\lnP{=}\prcondHLCPNm#7,\\#8&\prcondHLCPNm\lnP{=}\prcondHLCPNm#9.  \end{split}\end{equation} }
		\newcommand{\envHLineCECmDefNm}[9][*]{ \begin{equation}\begin{split} \ifx*#1 \text{[ERROR;need label name]} \else \label{#1} \fi #2&\prcondHLCPNm\lnP{:=}\prcondHLCPNm#3,\\#4&\prcondHLCPNm\lnP{:=}\prcondHLCPNm#5,\\#6&\prcondHLCPNm\lnP{:=}\prcondHLCPNm#7,\\#8&\prcondHLCPNm\lnP{:=}\prcondHLCPNm#9,  \end{split}\end{equation} }
		\newcommand{\envHLineCEDefNm}[9][*]{ \begin{equation}\begin{split} \ifx*#1 \text{[ERROR;need label name]} \else \label{#1} \fi #2&\prcondHLCPNm\lnP{:=}\prcondHLCPNm#3\\#4&\prcondHLCPNm\lnP{:=}\prcondHLCPNm#5\\#6&\prcondHLCPNm\lnP{:=}\prcondHLCPNm#7\\#8&\prcondHLCPNm\lnP{:=}\prcondHLCPNm#9  \end{split}\end{equation} }
		\newcommand{\envHLineCEDefNmPd}[9][*]{ \begin{equation}\begin{split} \ifx*#1 \text{[ERROR;need label name]} \else \label{#1} \fi #2&\prcondHLCPNm\lnP{:=}\prcondHLCPNm#3,\\#4&\prcondHLCPNm\lnP{:=}\prcondHLCPNm#5,\\#6&\prcondHLCPNm\lnP{:=}\prcondHLCPNm#7,\\#8&\prcondHLCPNm\lnP{:=}\prcondHLCPNm#9.  \end{split}\end{equation} }
	\newcommand{\envMLineTPd}[3][*]{ \ifx*#1 \begin{multline*} #2\lnP{=}#3.\end{multline*} \else \begin{multline} \label{#1} #2\lnP{=}#3.\end{multline} \fi }
	\newcommand{\envMLineTCm}[3][*]{ \ifx*#1 \begin{multline*} #2\lnP{=}#3,\end{multline*} \else \begin{multline} \label{#1} #2\lnP{=}#3,\end{multline} \fi }
	\newcommand{\envMLineTDefPd}[3][*]{ \ifx*#1 \begin{multline*} #2\lnP{:=}#3.\end{multline*} \else \begin{multline} \label{#1} #2\lnP{:=}#3.\end{multline} \fi }
	\newcommand{\envMLineTCmDef}[3][*]{ \ifx*#1 \begin{multline*} #2\lnP{:=}#3,\end{multline*} \else \begin{multline} \label{#1} #2\lnP{:=}#3,\end{multline} \fi }

	\newcommand{\envCaseTCm}[3][?]{\begin{cases} \glcondDis{#1}#2,\lcparaCase\\\glcondDis{#1}#3,\end{cases}}
	\newcommand{\envCaseTPd}[3][?]{\begin{cases} \glcondDis{#1}#2,\lcparaCase\\\glcondDis{#1}#3.\end{cases}}

	\DeclareFontEncoding{OT2}{}{}
	\DeclareFontSubstitution{OT2}{cmr}{m}{n}
	\DeclareFontFamily{OT2}{cmr}{\hyphenchar\font45}
	\DeclareFontShape{OT2}{cmr}{m}{n}{<5><6><7><8><9>gen*wncyr <10><10.95><12><14.4><17.28><20.74><24.88>wncyr10}{}
	\DeclareFontShape{OT2}{cmr}{b}{n}{<5><6><7><8><9>gen*wncyb<10><10.95><12><14.4><17.28><20.74><24.88>wncyb10}{}
	\DeclareMathAlphabet{\mathcyr}{OT2}{cmr}{m}{n}
	\DeclareMathAlphabet{\mathcyb}{OT2}{cmr}{b}{n}
	\SetMathAlphabet{\mathcyr}{bold}{OT2}{cmr}{b}{n}

	\newcommand{\shh}{*}
	\newcommand{\shs}{\mathcyr{sh}}	
			
	\newcommand{\shH}[1][\,]{#1\shh#1}
	\newcommand{\shS}[1][\;]{#1\shs#1}
	
		\newcommand{\sh}{\shs}
	\newcommand{\ztRLettHelp}{\bullet}
	\newcommand{\ztN}[1][?]{\zeta}				\newcommand{\ztO}[2][n]{\Fc[#1]{\zeta}{#2}}							
		\newcommand{\zt}{\ztO}
	\newcommand{\ztHN}[1][?]{\zeta_*}				\newcommand{\ztHO}[2][n]{\idFc[#1]{\zeta}{*}{#2}}					\newcommand{\ztHT}[3][n]{\idFc[#1]{\zeta}{*}{#2;#3}}		
		\newcommand{\ztH}{\ztHO}
	\newcommand{\ztRN}[1][?]{\zeta_\ztRLettHelp}						
					
	\newcommand{\ztSN}[1][?]{\zeta_\sh}			\newcommand{\ztSO}[2][n]{\idFc[#1]{\zeta}{\sh}{#2}}				\newcommand{\ztST}[3][n]{\idFc[#1]{\zeta}{\sh}{#2;#3}}	
		\newcommand{\ztS}{\ztSO}
	\newcommand{\ztsLettHelp}{\star}
	\newcommand{\ztsN}[1][?]{\zeta^\ztsLettHelp}				\newcommand{\ztsO}[2][n]{\pwFc[#1]{\zeta}{\ztsLettHelp}{#2}}							
		\newcommand{\zts}{\ztsO}
	\newcommand{\ztsHN}[1][?]{\zeta_*^\ztsLettHelp}				\newcommand{\ztsHO}[2][n]{\ipFc[#1]{\zeta}{*}{\ztsLettHelp}{#2}}					\newcommand{\ztsHT}[3][n]{\ipFc[#1]{\zeta}{*}{\ztsLettHelp}{#2;#3}}		
		\newcommand{\ztsH}{\ztsHO}
	\newcommand{\ztsRN}[1][?]{\zeta_\ztRLettHelp^\ztsLettHelp}					
					
	\newcommand{\ztsSN}[1][?]{\zeta_\sh^\ztsLettHelp}			\newcommand{\ztsSO}[2][n]{\ipFc[#1]{\zeta}{\sh}{\ztsLettHelp}{#2}}				\newcommand{\ztsST}[3][n]{\ipFc[#1]{\zeta}{\sh}{\ztsLettHelp}{#2;#3}}	
		\newcommand{\ztsS}{\ztsSO}

	\newcommand{\mpRG}[2][n]{\Fc[#1]{\rho}{#2}}						
				
	\newcommand{\mpsRG}[2][n]{\pwFc[#1]{\rho}{\star}{#2}}			
	\newcommand{\mpsRGi}[2][n]{\pwFc[#1]{{\rho^\star}}{ -1}{#2}}		

	\newcommand{\pintCO}[2][n]{\Fc[#1]{c}{#2}}	\newcommand{\pintCT}[3][n]{\idFc[#1]{c}{#2}{#3}}
	\newcommand{\pintsCO}[2][n]{\pwFc[#1]{c}{\star}{#2}}		\newcommand{\pintsCT}[3][n]{\ipFc[#1]{c}{#2}{\star}{#3}}
		\newcommand{\intCO}{\pintCO}
		\newcommand{\intsCO}{\pintsCO}
		\newcommand{\intCT}{\pintCT}
		\newcommand{\intsCT}{\pintsCT}
		\newcommand{\intC}{\intCO}
		\newcommand{\intsC}{\intsCO}
	\newcommand{\pletB}{\ztRLettHelp}
		\newcommand{\letB}{\pletB}
	\newcommand{\lettPztPIH}{H}
		\newcommand{\pztPIHTh}[4][n]{\idFc[#1]{\lettPztPIH}{\shH[]}{#2,#3;#4}}
		
		\newcommand{\ztPIHTh}{\pztPIHTh}
		\newcommand{\ztPIH}{\pztPIHTh}
		\newcommand{\pztPISTh}[4][n]{\idFc[#1]{\lettPztPIH}{\shS[]}{#2,#3;#4}}
		
		\newcommand{\ztPISTh}{\pztPISTh}
		\newcommand{\ztPIS}{\pztPISTh}
		\newcommand{\pztPIBTh}[4][n]{\idFc[#1]{\lettPztPIH}{\letB}{#2,#3;#4}}
		
		\newcommand{\ztPIBTh}{\pztPIBTh}
		
	\newcommand{\prefPropB}[2][?]{\refHL{Prop. }{#2}}
		\newcommand{\refPropB}{\prefPropB}
	
	\newcommand{\psetPI}[2][?]{\mathcal{P}_{#2}}
		\newcommand{\setPI}{\psetPI}
	\newcommand{\psetPTO}[2][n]{\Fc[#1]{\mathcal{P}}{#2}}	\newcommand{\psetPTT}[3][n]{\idFc[#1]{\mathcal{P}}{#2}{#3}}
		\newcommand{\setPTO}{\psetPTO}
		\newcommand{\setPTT}{\psetPTT}
		\newcommand{\setPT}{\setPTO}
	\newcommand{\psetS}[2][?]{[#2]}
		\newcommand{\setS}{\psetS}
		\newcommand{\pchKST}[3][n]{\idFc[#1]{\chi}{\shS[]}{#2,#3} }
		
		\newcommand{\chKST}{\pchKST}
		\newcommand{\chKS}{\chKST}
	\newcommand{\pfcET}[3][n]{\Fc[#1]{\eta}{#2;#3}}
		\newcommand{\fcET}{\pfcET}
	\newcommand{\pfcAO}[2][n]{\Fc[#1]{A}{#2}}
		
		\newcommand{\fcA}{\pfcAO}
	\newcommand{\pnmPT}[3][n]{\idFc[#1]{N}{#2}{#3}}
		\newcommand{\nmPT}{\pnmPT}
	\newcommand{\tpTdUn}[3][?]{ \dUn[#1]{ \tpT{#2}{#3} } }
		\newcommand{\tpdUn}{\tpTdUn}
	\newcommand{\pplBellT}[3][n]{\idFc[#1]{B}{#2}{#3}}
		\newcommand{\plBell}{\pplBellT}
	\newcommand{\pcoeBellT}[3][n]{\idFc[#1]{b}{#2}{#3}}
		\newcommand{\coeBell}{\pcoeBellT}
	\newcommand{\pnmSTfO}[2][n]{\Fc[#1]{c}{#2}}
		\newcommand{\nmSTf}{\pnmSTfO}

\allowdisplaybreaks[4]
	\title{\mainTitle}
	\author{\authorName
			\thanks{\organizationNameFst, \placeAddressFst}
			\mbox{}
			\thanks{\organizationNameScd, \departmentNameScd, \placeAddressScd}
		}
	\date{}

\begin{document}
\maketitle
\renewcommand{\thefootnote}{\fnsymbol{footnote}}
\footnote[0]{e-mail : \emailAddressFst}
\footnote[0]{MSC-class: \MSCname}
\renewcommand{\thefootnote}{\arabic{footnote}}\setcounter{footnote}{0}
\vPack[30]

\begin{abstract}
An identity involving symmetric sums of regularized multiple zeta-star values of harmonic type was proved
	by Hoffman. 
In this paper,
	we prove an identity of shuffle type.
We use Bell polynomials appearing  in the study of set partitions to prove the identity.
\end{abstract}
\section{Introduction and statement of results} \label{sectOne}
The multiple zeta value (MZV) and multiple zeta-star value (MZSV, or sometimes referred to as the non-strict MZV) 
	are real numbers defined by the nested series  
	\envHLine[1_PL_DefMZV_ser]
	{
		\zt{k_1, k_2, \ldots,k_r}
	}
	{
		\Sm{0<m_1 < m_2 < \cdots < m_r} \frc{m_1^{k_1} m_2^{k_2} \cdots m_r^{k_r}}{1}
	}
	and
	\envHLineCm[1_PL_DefMZSV_ser]
	{
		\zts{k_1, k_2, \ldots,k_r}
	}
	{
		\Sm{0<m_1 \leq m_2 \leq \cdots \leq m_r} \frc{m_1^{k_1} m_2^{k_2} \cdots m_r^{k_r}}{1}
	}
	respectively,
	where $k_i$ $(1 \leq i \leq r)$ are arbitrary positive integers with $k_r > 1$.
MZVs and MZSVs can also be given by integrals. 
The two expressions of series and integral yield
	two different products $\shH[]$ and $\shS[]$,  
	called \emph{harmonic} (or \emph{stuffle}) and \emph{shuffle},
	respectively (see, e.g., \cite{Hoffman97,IKZ06, KY18, Reutenauer93}).
These values have been actively studied since more than two decades,
	but 
	Euler \cite{Euler1776} already mentioned them in a special case, $r=2$.

We call a finite sequence $\bfK=(k_1, \ldots, k_r)$ of positive integers an \emph{index}.
MZVs and MZSVs are divergent if $k_r=1$,
	but recently,
	the theory of regularization has been established. 
(For details, see \cite{IKZ06} and \cite{KY18} for MZV and MZSV, respectively.) 
Four polynomials whose coefficients are $\setQ$-linear combinations of MZVs and MZSVs,
	which we denote by
	\envPLineCm
	{
		\ztHT{\bfK}{T},\quad \ztST{\bfK}{T}, \quad \ztsHT{\bfK}{T}, \quad\text{and}\quad \ztsST{\bfK}{T} 
	}
	are defined for any index $\bfK$ in the theory:
	$\ztHT{\bfK}{T}$ and $\ztST{\bfK}{T}$ are generalizations of $\zt{\bfK}$ involving products $\shH[]$ and $\shS[]$,
	respectively;
	$\ztsHT{\bfK}{T}$ and $\ztsST{\bfK}{T}$ are those of $\zts{\bfK}$.
A key idea of the generalizations is roughly 
	to regard the divergent value $\zt{1}=\frc[s]{1}{1}+\frc[s]{2}{1}+\cdots$ as the variable $T$ with keeping product rule. 
The regularized values $\ztH{\bfK}$, $\ztS{\bfK}$, $\ztsH{\bfK}$, and $\ztsS{\bfK}$ 
	are defined by their constant terms (e.g., $\ztH{\bfK}=\ztHT{\bfK}{0}$),
	and
	these values are of course equal to $\zt{\bfK}$ if $k_r>1$.
Fundamental theorems of regularization for MZVs and MZSVs were proved in \cite{IKZ06} and \cite{KY18}, respectively,
	which are stated as follows.
For any index $\bfK$, \hypertarget{1_PL_EqRegThms_hyper}{}
	\envPLineCm[1_PL_EqRegThms]
	{
		\mpRG{ \ztHT{\bfK}{T} }
	\lnP{=}
		\ztST{\bfK}{T}
	\qquad\text{and}\qquad
		\mpsRG{ \ztsHT{\bfK}{T} }
	\lnP{=}
		\ztsST{\bfK}{T}
	} 
	where $\mpRG[]{}$ and $\mpsRG[]{}$ are $\setR$-linear endomorphisms on $\setR[T]$ related with the gamma function $\fcG{u}$.
The detailed definition of $\mpsRG[]{}$ will be introduced in \refSec{sectTwo},
	which is necessary to prove our result, or \refThm{1_Thm1}.

In order to state \refThm{1_Thm1},
	we will recall Hoffman's identities involving symmetric sums of the polynomials $\ztHT{\bfK}{T}$ and $\ztsHT{\bfK}{T}$, 
	which are shown in \cite{Hoffman92,Hoffman15}.
Let $\setPI{r}$ be the set of partitions of the set $\SetO{1,\ldots,r}$. 
For any $\Pi = \bkB{P_1,\ldots,P_g} \in \setPI{r}$,
	we define integers $\intC{\Pi}=\intCT{r}{\Pi}$ and $\intsC{\Pi}=\intsCT{r}{\Pi}$ by
	\envPLineCm[1_PL_DefIntsC]
	{
		\intCT{r}{\Pi}
	\lnP{=}
		\mo^{r-g} \PdT{i=1}{g} (\vlNS{P_i}-1)! 
	\qquad\text{and}\qquad
		\intsCT{r}{\Pi}
	\lnP{=}	
		\PdT{i=1}{g} (\vlNS{P_i}-1)! 
	}
	respectively,
	where $\nmSet{P}$ is the number of the elements of a set $P$.
We also define 
	\envHLineCmDef[1_PL_DefZetaPartH]
	{
		\ztPIH{\bfK}{\Pi}{T}
	}
	{
		\PdT{i=1}{g} \fcET[G]{ \Sm{p\in P_i} k_p }{T}
	}
	where\footnote{We note that $\fcET{k}{T}=\ztHT{k}{T}=\ztST{k}{T}=\ztsHT{k}{T}=\ztsST{k}{T}$.}
	\envHLineDef
	{
		\fcET{k}{T}
	}
	{
		\envCaseTPd{
			\zt{k}
			&
			(k>1)
		}{
			T
			&
			(k=1)
		}
	} 
Let $\gpSym{r}$ denote the symmetric group of degree $r$.
Hoffman's identities are then
	\envHLine[1_PL_EqSymMZV]
	{
		\Sm{\sig\in\gpSym{r}} \ztHT{k_{\sig(1)}, \ldots, k_{\sig(r)}}{T}
	}
	{
		\Sm{\Pi\in\setPI{r}} \intC{\Pi} \ztPIHTh{\bfK}{\Pi}{T}
	}
	and
	\envHLinePd[1_PL_EqSymMZSV]
	{
		\Sm{\sig\in\gpSym{r}} \ztsHT{k_{\sig(1)}, \ldots, k_{\sig(r)}}{T}
	}
	{
		\Sm{\Pi\in\setPI{r}} \intsC{\Pi} \ztPIHTh{\bfK}{\Pi}{T}
	}
Recently, 
	a shuffle version of \refEq{1_PL_EqSymMZV} was proved in \cite{Machide17},
	which is obtained by replacing $\ztH[]{}$ and $\lettPztPIH_{\shH[]}$ with $\ztS[]{}$ and $\lettPztPIH_{\shS[]}$,
	respectively
	(see \refEq{1_Thm1_DefZetaPartS} for the definition of the function $\lettPztPIH_{\shS[]}$).

The main result of this paper is the shuffle version of \refEq{1_PL_EqSymMZSV}.

\begin{theorem}\label{1_Thm1}
For any index $\bfK$,
	we have
	\envHLineCm[1_Thm1_EqSymMZSV]
	{
		\Sm{\sig\in\gpSym{r}} \ztsST{k_{\sig(1)}, \ldots, k_{\sig(r)}}{T}
	}
	{
		\Sm{\Pi\in\setPI{r}} \intsC{\Pi} \ztPISTh{\bfK}{\Pi}{T} 
	}
	where $\ztPISTh{\bfK}{\Pi}{T}$ is similar to $\ztPIHTh{\bfK}{\Pi}{T}$,
	but the characteristic function 
	\envHLineDef[1_Thm1_DefCharFuncPIS]
	{
		\chKS{\bfK}{P_i}
	}
	{
		\envCaseTCm{
			0
			&
			\text{if $\nmSet{P_i}>1$, and $k_p=1$ for all $p\in P_i$}
		}{
			1
			&
			\text{otherwise}
		}
	}
	is added in each multiplicand;
	that is,
	\envHLineDefPd[1_Thm1_DefZetaPartS]
	{
		\ztPIS{\bfK}{\Pi}{T}
	}
	{
		\PdT{i=1}{g} \chKS{\bfK}{P_i} \fcET[G]{ \Sm{p\in P_i} k_p }{T}
	}
\end{theorem}

We give some examples of \refEq{1_Thm1_EqSymMZSV}.
The number of the terms of its right-hand side decreases
	as 
	the number of $k_i=1$ increases 
	because of \refEq{1_Thm1_DefCharFuncPIS}.

\begin{example}\label{1_Ex1}
Let $k$ and $l$ be integers at least $2$.
	\envHLineCm
	{
		\zts{1, k} + \ztsST{k, 1}{T}
	}
	{
		\zt{k}T + \zt{k+1}
	}
	\vPack[20]\envLineCm
	{
		\zts{1, k, l} + \zts{1, l, k} + \zts{k, 1, l} + \zts{l, 1, k} + \ztsST{k, l, 1}{T} + \ztsST{l, k, 1}{T} 
	}
	{
		\bkR{\zt{k}\zt{l}+\zt{k+l}}T + \zt{k}\zt{l+1} + \zt{l}\zt{k+1} + 2 \zt{k+l+1}
	}
	\vPack[20]\envHLinePd
	{
		2 \bkR{ \zts{1, 1, k} + \ztsST{1, k, 1}{T} + \ztsST{k, 1, 1}{T} }
	}
	{
		\zt{k}T^2 + 2\zt{k+1}T + 2\zt{k+2}
	}
\end{example}

In particular,
	we have a simple equation \refEq{1_Cor1_Eq}
	when the number of $k_i=1$ is $r-1$
	(or equivalently, there is just one $k_j$ that is greater than $1$):
	the right hand side are written in terms of only single zeta values. 
\begin{corollary}\label{1_Cor1}
For integers $k\geq 2$ and $r\geq1$,
	we have
	\envHLineCm[1_Cor1_Eq]
	{
		\SmT{i=0}{r-1} \ztsST{\bkB{1}^i, k, \bkB{1}^{r-1-i}}{T}
	}
	{
		\SmT{j=0}{r-1} \zt{k+r-1-j} \frc{j!}{T^{j}}	
	}	
	where $\bkB{1}^i$ means $i$ repetitions of $1$.
\end{corollary}

	
The method of the proof of \refThm{1_Thm1} is an improvement to that used  in \cite{Machide17}. 
We will use complete exponential Bell polynomials to show \refProp{2_Prop3},
	which are defined by 
	\envHLineDefPd[1_PL_DefBellPoly]
	{
		\plBell[a]{r}{x_1, \ldots, x_r}
	}
	{
		r! \tpSm{i_1,i_2,\ldots, i_r \geq0}{ 1\cdot i_1 + 2\cdot i_2 + \cdots + r\cdot i_r =r }
		\frc{i_1! i_2! \cdots i_r!}{ 1 } 
		\PdT{a=1}{r} \bkR[a]{ \frc{a!}{x_a} }^{i_a} 
	}
Bell polynomials \cite{Bell27} appear in the study of set partitions at first.
Currently it is known that
	they have many relations to combinatorial numbers and applications to other areas (see, e.g., \cite{Comtet74, Roman80}).
We will mention an identity involving $\ztsHT{1, 1, \ldots, 1}{T}$ in \refRem{2_Rem1},
	which appears a variation of the identity 
	\envMPd
	{
		r!
	}
	{
		\plBell[a]{r}{0!, 1!, \ldots, (r-1)!}
	}

This paper is organized as follows.
We prepare some propositions in \refSec{sectTwo},
	and
	prove \refThm{1_Thm1} and \refCor{1_Cor1} in \refSec{sectThree}.

\section{Propositions} \label{sectTwo}
In this section,
	we introduce \refProp[s]{2_Prop1}, \ref{2_Prop2}, and \ref{2_Prop3}
	that will be used to prove \refThm{1_Thm1}.
We will omit the proofs of \refProp[s]{2_Prop1} and \ref{2_Prop2}
	because these are almost the same as \refLemA[s]{4.7} and {4.8} in \cite{Machide17},
	respectively,
	where some notation and terminology are modified.

Let $\setS{r}$ denotes the set $\bkB{1, \ldots, r}$,
	and 
	let $A$ and $B$ be its subsets.
We denote by $\setPT{A}$ the set of partitions of $A$ (i.e., $\setPT{\setS{r}}=\setPI{r}$),
	and
	we define a subset $\setPTT{B}{A}$ in $\setPT{A}$ by
	\envHLineDefPd
	{
		\setPTT{B}{A}
	}
	{
		\SetT{ \Pi = \Set{P_1,\ldots,P_m} \in \setPT{A} }{\text{$P_i \not\subset B$ for all $i$}}
	}
For example,
	if $(r, A, B) = (4, \Set{1,2,3}, \Set{3,4})$,	
	then
	\envPLineCm
	{
		\setPT{A}
	\lnP{=}
		\Set[b]{ 1|2|3, 12|3, 13|2, 23|1, 123}
	\quad\text{and}\quad
		\setPTT{B}{A}
	\lnP{=}
		\Set[b]{ 13|2, 23|1, 123 }
	}
	where $a_1 \cdots a_p | b_1 \cdots b_q | \cdots$ means a partition
	such as
	\envMPd
	{
		12|3
	}
	{
		\Set{ \Set{1,2}, \Set{3} }
	}

Let $\Xi=\Set{P_1,\ldots,P_g} \in \setPT{A}$,
	and let $s=\nmSet{A}$.
We will define a partition $\sig_A(\Xi)$ in $\setPI{s}$, 
	as follows.
Let $a_1<\cdots<a_s$ be the increasing sequence of integers such that 
	\envHLineCm
	{
		A
	}{
		\Set{a_1,\ldots,a_s}
	}
	and
	let $\sig_A$ be the permutation of $\gpSym{r}$ that is uniquely determined by  
	\envPLine
	{
		\sig_A^{-1}(i)
	\lnP{=}
		a_i	\;(i=1,\ldots,s)
	\qquad\text{and}\qquad
		\sig_A^{-1}(s+1) < \cdots < \sig_A^{-1}(r)
	;}
	by the definition,
	\envOTLineThPd
	{
		\sig_A(A)
	}
	{
		\Set{\sig_A(a_{1}), \ldots, \sig_A(a_{s})}
	}
	{
		\setS{s}
	}
We then define 
	\envHLineDefPd
	{
		\sig_A(\Xi)
	}
	{
		\Set{\sig_A(P_1),\ldots,\sig_A(P_g)}
	\lnP{\in}
		\setPI{s}
	}
For convenience, 
	$\sig_A(\Xi) = \phi$ if $A=\Xi=\phi$.  

The propositions are as follows.

\begin{proposition}[{see \cite[Lemma 4.7]{Machide17}}]\label{2_Prop1}
For any subset $B \subsetneqq \setS{r}$,
	we have
	\envHLineCm[2_Prop1_Eq1]
	{
		\dUn{A \subset B} \SetT{ \Xi \sqcup \Del }{ (\Xi, \Del) \in \setPT{A} \times \setPTT{B}{ \setS{r} \setminus A } }
	}
	{
		\setPI{r}
	}
	where 
	$\sqcup$ denotes the disjoint union,
	and
	$\dUn{A \subset B}$ ranges over all subset in $B$ which include $\phi$.
\end{proposition}

\begin{proposition}[{see \cite[Lemma 4.8]{Machide17}}]\label{2_Prop2}
Let $A$ and $B$ be subsets such that $A \subset B \subsetneqq \setS{r}$,
	and 
	let $(\Xi, \Del)$ be in $\setPT{A} \times \setPTT{B}{ \setS{r} \setminus A }$.
Let the symbol $\letB$ mean either $\shH[]$ or $\shS[]$. 
\mbox{}\\{\bf(i)}
We define $\intsC{\phi}=1$.
We have
	\envHLinePd[2_Prop2_Eq1]
	{
		\intsC{\Xi \cup \Del}
	}
	{
		\intsC{\Xi} \intsC{\Del} 
	}
{\bf(ii)}
We define $\ztPIBTh{\phi}{\phi}{T}=1$ and $\mathbf{1}_s=( \usbTx{s}{1,\ldots,1})$.
Suppose that 
	$\bfK=(k_1, \ldots, k_r)$ is an index satisfying  
	\envM
	{
		B
	}
	{ 
		\SetT{ a \in \setS{r} }{ k_a = 1}
	}
	(and $\bfK\neq\mathbf{1}_r$). 
Then we have
	\envHLineCm[2_Prop2_Eq2]
	{
		\ztPIBTh{\bfK}{\Xi \cup \Del}{T}
	}
	{
		\bkR[G]{ \PdT{i=1}{h} \zt{ k_{Q_i} } }  \ztPIBTh{ \mathbf{1}_{\nmSet{A}} }{\sig_A(\Xi)}{T}
	}
	where
	$Q_1, \ldots, Q_h$ are the blocks of $\Del$ (i.e., $\Del = \SetO{ Q_1, \ldots, Q_h }$),
	and
	\envHLinePd
	{
		k_{Q_i}
	} 
	{
		\Sm{q\in Q_i} k_q
		\qquad
		(i=1,\ldots,h)
	}
Note that
	$k_{Q_i}>1$ and $\zt{ k_{Q_i} }$ is not infinity for any $i$, 
	because $\Del \in \setPTT{B}{ \setS{r} \setminus A }$ and $Q_i \not\subset B$.
\end{proposition}

\begin{proposition}[{see \cite[Lemma 4.9]{Machide17}}]\label{2_Prop3}
For any positive integer $r$,
	we have
	\envHLineCFNmePd
	{\label{2_Prop3_Eq1}
		\Sm{ \Pi \in \setPI{r} } \intsC{\Pi} \ztPIHTh{ \mathbf{1}_r }{\Pi}{T}
	}
	{
		\mpsRGi{T^r}
	}
	{\label{2_Prop3_Eq2}
		\Sm{ \Pi \in \setPI{r} } \intsC{\Pi} \ztPISTh{ \mathbf{1}_r }{\Pi}{T}
	}
	{
		T^r
	}
\end{proposition}

The condition $B \neq \setS{r}$ in the first two propositions is necessary for taking an element in $\setPTT{B}{ \setS{r} \setminus A }$
	(see \cite[Remark 4.6]{Machide17} for details).

To prove \refProp{2_Prop3},
	we require \refLem{2_Lem1} 
	that is the \emph{star}-version of \cite[\refLemA{4.10}]{Machide17} in terms of Bell polynomials.

\begin{lemma}\label{2_Lem1}		
For any positive integer $r$,
	we have
	\envHLinePd[2_Lem1_Eq]
	{
		\Sm{\Pi\in\setPI{r}} \intsCT{}{\Pi} \ztPIHTh{\mathbf{1}_r}{\Pi}{T}
	}
	{
		\plBell[a]{r}{0!\fcET{1}{T}, 1!\fcET{2}{T}, \ldots, (r-1)!\fcET{r}{T}}
	}
\end{lemma}

We will now prove \refProp{2_Prop3},
	and then prove \refLem{2_Lem1}.

\envProof[\refProp{2_Prop3}]{
We first recall the definition of $\mpsRG[]{}$,
	which is an $\setR$-linear endomorphism on $\setR[T]$ determined by the equality 
	\envHLine[2_Prop3Pr_DefMpsRG]
	{
		\mpsRG{\exN{Tt}}
	}
	{
		\fcA{-t}^{-1} \exN{Tt}
	} 
	in the formal power series algebra $\setR[T][[t]]$ on which $\mpsRG[]{}$ acts coefficientwise (see \cite[\refSecA{4}]{KY18}),
	where 
	\envHLinePd
	{
		\fcA{t}
	}
	{
		\exx[G]{\SmT{m=2}{\infty} \frc{m}{\mo^m\zt{m}} t^m }
	}
Note that 
	$\fcA{t}=\exN{\gam t}\fcG{1+t}$,
	where $\gam$ is Euler's constant.
We can see from \refEq{2_Prop3Pr_DefMpsRG} 
	that 
	the inverse endomorphism $\mpsRGi[]{}$ exists 
	and 
	it satisfies
	\envHLineFPd[2_Prop3Pr_Eq1]
	{
		\mpsRGi{\exN{Tt}}
	}
	{
		\fcA{-t} \exN{Tt} 
	}
	{
		\exx[a]{ Tt + \SmT{m=2}{\infty} \frc{ m }{ \zt{m} } t^m }
	}
	{
		\exx[a]{ \SmT{m=1}{\infty} \frc{ m }{ \fcET{m}{T} } t^m }
	}
The exponential partial Bell polynomials can be defined by use of the generating function (see \cite[Chapter 3]{Comtet74}):
	\envHLinePd[2_Prop3Pr_Eq2]
	{
		\exx[a]{ \SmT{m=1}{\infty} x_m\frc{m!}{t^m} }
	}
	{
		\SmT{r=0}{\infty} \plBell{r}{x_1, \ldots, x_r} \frc{r!}{t^r}
	}
Combining \refEq{2_Prop3Pr_Eq1} and \refEq{2_Prop3Pr_Eq2} with $x_m= (m-1)! \fcET{m}{T}$,
	we obtain
	\envHLineCm
	{
		\mpsRGi{\exN{Tt}}
	}
	{
		\SmT{r=0}{\infty} \plBell[a]{r}{0!\fcET{1}{T}, 1!\fcET{2}{T}, \ldots, (r-1)!\fcET{r}{T}} \frc{r!}{t^r} 
	}
	which,
	together with \refEq{2_Lem1_Eq},
	gives
	\envHLinePd
	{
		\mpsRGi{\exN{Tt}}	
	}
	{
		\SmT{r=0}{\infty} \frc{r!}{t^r} \Sm{ \Pi \in \setPI{r} } \intsC{\Pi} \ztPIHTh{ \mathbf{1}_r }{\Pi}{T}
	}
Identity \refEq{2_Prop3_Eq1} follows from comparing the coefficients of $t^r$ on both sides of this equation. 

We will give a proof of \refEq{2_Prop3_Eq2},
	which is a modification to that of \refEqA{4.58} 
	in \cite[Lemma 4.9]{Machide17}.\footnote{There is a misprint in the proof of \cite[\refEqA{4.58}]{Machide17}: 
			$\SetO{ \usbTx{n}{\SetO{1}, \ldots, \SetO{1} }}$ should be $\SetO{ \SetO{1}, \ldots, \SetO{n} }$.}
Let $\Lambda=\Lambda_r$ be the partition in $\setPI{r}$ defined by
	\envHLineDefPd
	{
		\Lambda_r
	}
	{
		1|2|\cdots|r
	\lnP{=}
		\SetO{ \SetO{1}, \SetO{2}, \ldots, \SetO{r} }
	} 
We see from \refEq{1_Thm1_DefCharFuncPIS} and \refEq{1_Thm1_DefZetaPartS} that
	\envM
	{
		\ztPISTh{ \mathbf{1}_r }{ \Pi }{T}
	}
	{
		0
	}
	for any $\Pi \in \setPI{r}$ with $\Pi \neq \Lambda$,
	and so 
	\envHLinePd
	{
		\Sm{ \Pi \in \setPI{r} } \intsC{\Pi} \ztPISTh{ \mathbf{1}_r }{\Pi}{T}
	}
	{
		\intsC{\Lambda} \ztPISTh{ \mathbf{1}_r }{ \Lambda }{T}
	}
Since
	\envPLineCm
	{
		\intsC{\Lambda}
	\lnP{=}
		\PdT{i=1}{r} 0! 
	\lnP{=}
		1
	\qquad\text{and}\qquad
		\ztPISTh{ \mathbf{1}_r }{ \Lambda }{T}
	\lnP{=}
		\PdT{i=1}{r} \fcET{ 1}{T}
	\lnP{=}
		T^r
	}
	we obtain \refEq{2_Prop3_Eq2}.
}

We will need the partial exponential Bell polynomials to prove \refLem{2_Lem1},
	which we denote by $\plBell{r,k}{x_1, \ldots, x_{r-k+1}}$ for integers $r$ and $k$ with $1 \leq k \leq r$.
Complete and partial Bell polynomials have the relations
	\envHLinePd[2_PL_Eq1BellPolys]
	{
		\plBell{r}{x_1, \ldots, x_{r}}
	}
	{
		\SmT{k=1}{r} \plBell{r,k}{x_1, \ldots, x_{r-k+1}}
	}
Let $\coeBell{r,k}{i_1, \ldots, i_{r-k+1}}$ be the coefficients of $\plBell{r,k}{x_1, \ldots, x_{r-k+1}}$ such that
	\envHLinePd
	{
		\plBell{r,k}{x_1, \ldots, x_{r-k+1}}
	}
	{
		\tpSm{i_1,\ldots,i_{r-k+1}\geq0}{ i_1+i_2+\cdots+i_{r-k+1}=k \atop 1\cdot i_1 + 2\cdot i_2 + \cdots + (r-k+1)\cdot i_{r-k+1} =r  } 
		\coeBell{r,k}{i_1, \ldots, i_{r-k+1}} \PdT{a=1}{r-k+1} x_a^{i_a}
	}
From combinatorial consideration (see, e.g., \cite[Chapter 3]{Comtet74}),
	we know that
	$\coeBell{r,k}{i_1, \ldots, i_{r-k+1}}$  is the number of partitions with total $k$ blocks in $\setPI{r}$ 
	which 
	consist of $i_{a}$ blocks of size $a$ for $a \in \setS{r-k+1}$.	
For instance,
	\envPLineCm
	{
		\coeBell{4,2}{1,0,1}	\lnP{=}	4
	,\quad
		\coeBell{4,2}{0,2,0}	\lnP{=}	3
	,\quad\text{and}\quad
		\plBell{4,2}{x_1, x_2, x_3}	\lnP{=}	4x_1x_3 + 3x_2^2
	}
	since
	we have $4$ partitions with blocks of size $1$ and $3$ 
	and
	$3$ partitions with $2$ blocks of size $2$ 
	when a set with $4$ elements is divided into $2$ blocks.
We note that
	$\plBell[]{r,k}=\plBell{r,k}{1, \ldots, 1}$ are Stirling numbers of the second kind,
	that is,
	count the number of ways to partition a set of $r$ elements into $k$ nonempty subsets. 

\envProof[\refLem{2_Lem1}]{
For any partition $\Pi=\Set{P_1,\ldots,P_{g}}$ in $\setPI{r}$ and integer $a$ in $\setS{r}$,
	let $\nmPT{a}{\Pi}$ be the number of the blocks whose cardinalities equal $a$,
	i.e.,
	\envHLineDefPd
	{
		\nmPT{a}{\Pi}
	}
	{
		\nmSet{ \SetT{j\in\setS{g} }{ \nmSet{P_j}=a } }
	}
We see from the definition that
	\envHLineCFLaaCm
	{ 
		g
	}
	{
		\nmPT{1}{\Pi}+\cdots+\nmPT{r}{\Pi}
	}
	{
		r
	}
	{
		1\cdot \nmPT{1}{\Pi} + \cdots + r\cdot \nmPT{r}{\Pi}
	}
	so that 
	\envOTLinePd
	{
		\intsC{\Pi} \ztPIHTh{\mathbf{1}_r}{\Pi}{T}
	}
	{
		\PdT{i=1}{g} (\nmSet{P_i}-1)!  \fcET{\nmSet{P_i}}{T} 
	}
	{
		\PdT{a=1}{r} \bkR{ (a-1)! \fcET{ a }{T} }^{\nmPT{a}{\Pi}}
	}
It follows from the combinatorial meaning of $\coeBell{r,k}{i_1, \ldots, i_{r-k+1}}$ that
	\envHLine
	{
		\tpSm{\Pi\in\setPI{r}}{\nmPT{a}{\Pi}=i_a (\forall a)} 1
	}
	{
		\coeBell{r,k}{i_1, \ldots, i_{r-k+1}}
	}
	for non-negative integers $i_1, \ldots, i_{r-k+1}$ with $\SmT[l]{a=1}{r-k+1} i_a=k$ and $\SmT[l]{a=1}{r-k+1} a\cdot i_a=r$,
	and so
	\envLineFPd[2_Lem1Pr_Eq]
	{
		\tpSm{i_1,i_2,\ldots,i_r \geq 0}{1\cdot i_1 + 2\cdot i_2 + \cdots + r\cdot i_r=r} \tpSm{\Pi\in\setPI{r}}{\nmPT{a}{\Pi}=i_a (\forall a)}
		\intsC{\Pi} \ztPIHTh{ \mathbf{1}_r }{\Pi}{T}
	}
	{
		\tpSm{i_1,i_2,\ldots,i_r \geq 0}{1\cdot i_1+2\cdot i_2 + \cdots + r\cdot i_r=r} 
		\bkR[a]{ \PdT{a=1}{r} \bkR{ (a-1)! \fcET{a}{T} }^{i_a} }
		\tpSm{\Pi\in\setPI{r}}{\nmPT{a}{\Pi}=i_a (\forall a)} 1
	}
	{
		\SmT{k=1}{r}
		\tpSm{i_1,i_2,\ldots,i_{r-k+1} \geq 0}{i_1+i_2+\cdots+i_{r-k+1}=k \atop 1\cdot i_1+2\cdot i_2 + \cdots + (r-k+1)\cdot i_{r-k+1}=r} 
		\bkR[a]{ \PdT{a=1}{r-k+1} \bkR{ (a-1)! \fcET{a}{T} }^{i_a} }
		\coeBell{r,k}{i_1, \ldots, i_{r-k+1}}
	}
	{
		\SmT{k=1}{r} \plBell[a]{r,k}{0!\fcET{1}{T}, 1!\fcET{2}{T}, \ldots, (r-k)!\fcET{r-k+1}{T}}
	}
We thus obtain \refEq{2_Lem1_Eq},
	since the first line of \refEq{2_Lem1Pr_Eq} is equal to the left-hand side of \refEq{2_Lem1_Eq} because of \cite[\refEqA{4.74}]{Machide17},
	or
	\envHLineCm
	{
		\setPI{r}
	}
	{
		\tpdUn{i_1,i_2,\ldots,i_r \geq 0}{1\cdot i_1 + 2\cdot i_2 + \cdots + r\cdot i_r=r} 
		\SetT{ \Pi\in\setPI{r} }{ \nmPT{a}{\Pi}=i_a\,(a \in \setS{r}) }
	}
	and 
	since the last line of \refEq{2_Lem1Pr_Eq} is equal to the right-hand side of \refEq{2_Lem1_Eq} because of \refEq{2_PL_Eq1BellPolys}.
}
\begin{remark}\label{2_Rem1} 
Bell polynomials are related to many combinatorial numbers.
It may be worth noting the relation to the unsigned Stirling numbers of the first kind,
	which can be expressed as
	\envHLinePd{
		\nmSTf{r,k}
	}
	{
		\plBell{r,k}{0!, \ldots, (r-k)!}
	}
The unsigned Stirling numbers are defined as coefficients of the rising factorial,
	that is,
	\envHLinePd[2_Rem1_Eq0]
	{
		x(x+1)\cdots(x+r-1)
	}
	{
		\SmT[l]{k=1}{r} \plBell{r,k}{0!, \ldots, (r-k)!} x^k
	}
Substituting $x=1$ in this equation and combining it with \refEq{2_PL_Eq1BellPolys},
	we thus obtain
	\envHLinePd[2_Rem1_Eq1]
	{
		r!
	}
	{
		\plBell[a]{r}{0!, 1!, \ldots, (r-1)!}
	}
We also have
	\envHLineCm[2_Rem1_Eq2]
	{
		r! \ztsHT{\mathbf{1}_r}{T}
	}
	{
		\plBell[a]{r}{0!\fcET{1}{T}, 1!\fcET{2}{T}, \ldots, (r-1)!\fcET{r}{T}}
	}
	which follows from Hoffman's identity \refEq{1_PL_EqSymMZSV} with $\bfK=\mathbf{1}_r$ and \refEq{2_Lem1_Eq}.
Equation \refEq{2_Rem1_Eq2} is a variation of \refEq{2_Rem1_Eq1}
	in the sense that 
	we can obtain \refEq{2_Rem1_Eq2} from  \refEq{2_Rem1_Eq1} 
	by replacing $r!$ in the left-hand side with $r!\ztsST{\mathbf{1}_r}{T}$ and $j!$ in the right-hand side with $j!\fcET{j+1}{T}$.


\end{remark}	

\section{Proof} \label{sectThree}
We will need \refEq{3_Prop1_Eq} to prove \refEq{1_Thm1_EqSymMZSV},
	which is the star version of \cite[\refEqA{4.51}]{Machide17}.  

\begin{proposition}\label{3_Prop1}
For any index $\bfK$,
	\envHLinePd[3_Prop1_Eq]
	{
		\mpsRG[G]{ \Sm{\Pi\in\setPI{r}} \intsC{\Pi} \ztPIHTh{\bfK}{\Pi}{T}  }
	}
	{
		\Sm{\Pi\in\setPI{r}} \intsC{\Pi} \ztPISTh{\bfK}{\Pi}{T} 
	}
\end{proposition}

We can prove \refEq{3_Prop1_Eq} in a quite similar way 
	as we see below.
	

\envProof[\refProp{3_Prop1}]{
Let $B=\SetT{ j \in \setS{r} }{ k_j=1 }\subset\setS{r}$.
We suppose that $B=\setS{r}$.
Then, 
	$\bfK=\mathbf{1}_r$, 
	and so
	we can see from \refProp{2_Prop3} that
	\envHLineFCmPte[3_Prop1Pr_Eq1]
	{
		\mpsRG[G]{ \Sm{\Pi\in\setPI{r}} \intsC{\Pi} \ztPIHTh{\bfK}{\Pi}{T}  }
	}{\osTx{=}{\refEq{2_Prop3_Eq1}}}
	{
		\mpsRG{ \mpsRGi{T^r} }
	}{=}
	{
		T^r
	}{\osTx{=}{\refEq{2_Prop3_Eq2}}}
	{
		\Sm{ \Pi \in \setPI{r} } \intsC{\Pi} \ztPISTh{ \bfK }{\Pi}{T}
	}
	which proves \refEq{3_Prop1_Eq} for this case.

We suppose that $B\neq\setS{r}$.
Let $A$ be a subset in $B$.
Then we have
	\envHLineCm[3_Prop1Pr_Eq2]
	{
		\SetT{ \sig_A(\Xi) }{ \Xi \in \setPT{A} }
	}
	{
		\SetT{ \Xi' }{ \Xi' \in \setPI{\nmSet{A}} }
	}
	because the restriction of $\sig_A$ to $A$ is a bijection from $A$ to $\setS{\nmSet{A}}$.
From the definition \refEq{1_PL_DefIntsC}
	we easily see that
	\envMPd
	{
		\intsC{\Xi}
	}
	{
		\intsC{\sig_A(\Xi)}
	}
Thus,
	\envLineCmPt{\osTx{=}{(\refPropB{2_Prop1})}}
	{
		\Sm{\Pi\in\setPI{r}} \intsC{\Pi} \ztPIHTh{\bfK}{\Pi}{T}
	}
	{
		\Sm{A\subset B} \Sm{ {\Xi \in \setPT{A} } \atop { \Delta\in \setPTT{B}{ \setS{r}\setminus A } } } 
		\intsC{\Xi \cup \Delta} \ztPIHTh{\bfK}{\Xi \cup \Delta}{T}
	\lnAHP{\osTx{=}{(\refPropB{2_Prop2})}}
		\Sm{A\subset B} \Sm{ \Delta \in \setPTT{B}{ \setS{r}\setminus A}  }  
		\intsC{\Delta} \bkR[a]{ \PdT{i=1}{h}  \zt{k_{Q_i}} }
		\Sm{ \Xi \in \setPT{A} } \intsC{\Xi} \ztPIHTh{\mathbf{1}_{\nmSet{A}}}{\sig_A(\Xi)}{T}
	\lnAHP{\osTx{=}{(\ref{3_Prop1Pr_Eq2})}}
		\Sm{A\subset B} \Sm{ \Delta \in \setPTT{B}{ \setS{r}\setminus A }  }  
		\intsC{\Delta} \bkR[a]{ \PdT{i=1}{h}  \zt{k_{Q_i}} }
		\Sm{ \Xi' \in \setPI{\nmSet{A}} } \intsC{\Xi'} \ztPIHTh{\mathbf{1}_{\nmSet{A}}}{\Xi'}{T}
	\lnAHP{\osTx{=}{\refEq{2_Prop3_Eq1}}}
		\Sm{A\subset B} \Sm{ \Delta \in \setPTT{B}{ \setS{r}\setminus A }  }  
		\intsC{\Delta} \bkR[a]{ \PdT{i=1}{h}  \zt{k_{Q_i}} } 
		\mpsRGi{T^{\nmSet{A}}}
	}
	where $Q_1,\ldots,Q_h$ mean the blocks of $\Delta$. 
Therefore,
	\envLineThPd[3_Prop1Pr_Eq3]
	{
		\mpsRG[G]{ \Sm{\Pi\in\setPI{r}} \intsC{\Pi} \ztPIHTh{\bfK}{\Pi}{T}  }
	}
	{
		\Sm{A\subset B} \Sm{ \Delta \in \setPTT{B}{ \setS{r}\setminus A }  }  
		\intsC{\Delta} \bkR[a]{ \PdT{i=1}{h}  \zt{k_{Q_i}} } 
		\mpsRG{ \mpsRGi{T^{\nmSet{A}}} }
	}
	{
		\Sm{A\subset B} \Sm{ \Delta \in \setPTT{B}{ \setS{r}\setminus A }  }  
		\intsC{\Delta} \bkR[a]{ \PdT{i=1}{h}  \zt{k_{Q_i}} } 
		T^{\nmSet{A}}
	}
By using \refProp{2_Prop1}, \refProp{2_Prop2}, and \refEq{3_Prop1Pr_Eq2}, 
	and 
	by using \refEq{2_Prop3_Eq2} instead of \refEq{2_Prop3_Eq1},
	we can similarly prove
	\envHLinePd[3_Prop1Pr_Eq4]
	{
		\Sm{\Pi\in\setPI{r}} \intsC{\Pi} \ztPISTh{\bfK}{\Pi}{T}
	}
	{
		\Sm{A\subset B} \Sm{ \Delta \in \setPTT{B}{ \setS{r}\setminus A }  }  
		\intsC{\Delta} \bkR[a]{ \PdT{i=1}{h}  \zt{k_{Q_i}} } 
		T^{\nmSet{A}}
	}
Equating \refEq{3_Prop1Pr_Eq3} and \refEq{3_Prop1Pr_Eq4},
	we obtain \refEq{3_Prop1_Eq} for $B\neq\setS{r}$,
	and complete the proof.
}

We are now in a position to prove \refThm{1_Thm1}. 

\envProof[\refThm{1_Thm1}]{
Recall the \emph{star-regularization theorem}, the second identity of \refEq{1_PL_EqRegThms},
	which we tag at (\bfTx{s-reg}) here.
We obtain
	\envHLineCmPt{\osTx{=}{(\hyperlink{1_PL_EqRegThms_hyper}{\text{\bfTx{s-reg}}})}}
	{
		\Sm{\sig\in\gpSym{r}} \ztsST{k_{\sig(1)}, \ldots, k_{\sig(r)}}{T}
	}
	{
		\mpsRG[G]{ \Sm{\sig\in\gpSym{r}} \ztsHT{k_{\sig(1)}, \ldots, k_{\sig(r)}}{T} }
	\lnAHP{\osTx{=}{(\ref{1_PL_EqSymMZSV})}}
		\mpsRG[G]{ \Sm{\Pi\in\setPI{r}} \intsC{\Pi} \ztPIHTh{\bfK}{\Pi}{T} }
	\lnAHP{\osTx{=}{(\ref{3_Prop1_Eq})}}
		\Sm{\Pi\in\setPI{r}} \intsC{\Pi} \ztPISTh{\bfK}{\Pi}{T} 
	}
	which proves \refEq{1_Thm1_EqSymMZSV}.
}

Finally,
	we deduce \refCor{1_Cor1} from \refThm{1_Thm1}.

\envProof[\refCor{1_Cor1}]{
Let $\bfK=(k_1, \ldots, k_r)$ be the index $(k,1,\ldots,1)$,
	and
	let $\Pi = \SetO{P_1, \ldots, P_m}$ denote a partiton in $\setPI{r}$.
Assume that $1\in P_1$ through this proof,
	which dose not loss the generality.
We see from \refEq{1_Thm1_DefCharFuncPIS} that
	\envHLine
	{
		\ztPIS{\bfK}{\Pi}{T}
	}
	{
		\envCaseTPd{
			0
			&
			\text{if $\nmSet{P_i}>1$ for some $i \geq 2$}
		}{
			\zt{k+ \nmSet{P_1}-1} T^{m-1}
			&
			\text{otherwise}
		}
	}	
Since 
	\envOTLineTh
	{
		m-1
	}
	{
		\SmT{i=2}{m} \nmSet{P_i} 
	}
	{
		r - \nmSet{P_1} 
	}
	if $\nmSet{P_i}=1$ for all $i \geq 2$,
	it follows from \refEq{1_Thm1_EqSymMZSV} that
	\envLineThCm[3_1_Cor1Pr_Eq1]
	{
		(r-1)! \SmT{i=0}{r-1} \ztsST{\bkB{1}^i, k, \bkB{1}^{r-1-i}}{T}
	}
	{
		\SmT{j=1}{r} \Sm{\Pi \in \mathcal{X}_{j}} (j-1)! \zt{k+ j-1} T^{r-j}
	}
	{
		\SmT{j=1}{r} (j-1)! \zt{k+j-1} T^{r-j} \Sm{\Pi \in \mathcal{X}_{j}} 1
	}
	where
	\envHLinePd
	{
		\mathcal{X}_j
	}
	{
		\SetT{ \SetO{P_1, P_2, \ldots, P_{r+1-j}}\in\setPI{r} }{ \nmSet{P_1}=j, \nmSet{P_2}=\cdots=\nmSet{P_{r+1-j}}=1}
	}
Noting the assumption $1\in P_1$,
	we have 
	\envHLineCm
	{
		\mathcal{X}_j
	}
	{
		\SetT[b]{ \SetO{\setS{r}\setminus\SetO{a_2, \ldots, a_{r+1-j}}, \SetO{a_2}, \ldots, \SetO{a_{r+1-j}}} }{ 2 \leq a_2 < \cdots < a_{r+1-j} \leq r }
	}
	where 
	$\setS{r}\setminus\SetO{a_2, \ldots, a_{r+1-j}}$ corresponds to $P_1$
	and
	$\SetO{a_i}$ ($2 \leq i \leq r$) correspond to $P_i$.
Thus
	$\nmSet{ \mathcal{X}_j }$ is the number of $(r-j)$-combinations of $\SetO{2, \ldots, r}$,
	or
	\envHLinePd[3_1_Cor1Pr_Eq2]
	{
		\Sm{\Pi \in \mathcal{X}_{j}} 1
	}
	{
		\bnm{r-j}{r-1}
	}
Combining \refEq{3_1_Cor1Pr_Eq1} and \refEq{3_1_Cor1Pr_Eq2},
	we obtain
	\envHLinePd
	{
		\SmT{i=0}{r-1} \ztsST{\bkB{1}^i, k, \bkB{1}^{r-1-i}}{T}
	}
	{
		\SmT{j=1}{r} \frc{(r-j)!}{1} \zt{k+j-1} T^{r-j}
	}
Replacing $j$ with $r-j$ in the right-hand side of this equation gives \refEq{1_Cor1_Eq}.
}


\section*{Acknowledgements}
This work was supported by JST ERATO Grant Number JPMJER1201, Japan.



\end{document}